**Some sequences and number triangles which are related to**

**Narayana polynomials and to q-Narayana polynomials for q=-1.**

Johann Cigler

**Abstract**

As is well known the Catalan numbers $C_n$ can be interpreted as the number of Dyck paths of length $2n.$ The Narayana polynomials and the q-Narayana polynomials for q=-1 can analogously be interpreted as certain weights of these paths. In the present note we try to get some information about the corresponding weights of bounded Dyck paths. This also leads to some number sequences and analogues of Pascal's triangle which previously have occurred in other contexts. We consider some examples and obtain some results and conjectures.

**1. Introduction and background material**

Let

(1) $$C_n(t)=\sum_{k=0}^{n-1}\binom{n}{k}\binom{n-1}{k}\frac{t^k}{k+1}$$

be the Narayana polynomials

$$\left(C_n(t)\right)_{n\geq 1}=\left(1,1+t,1+3t+t^2,1+6t+6t^2+t^3,\cdots\right)$$

and

(2) $$C_n\left(t;q\right)=\sum_{k=0}^{n-1}q^{k(k+1)}\begin{bmatrix}n\\k\end{bmatrix}_q\begin{bmatrix}n-1\\k\end{bmatrix}_q\frac{t^k}{[k+1]_q}$$

$$\left(C_n(t;q)\right)_{n\geq 1}=\left(1,1+q^2t,1+q^2\left(1+q+q^2\right)t+q^6t^2,1+q^2\left(1+q^2\right)\left(1+q+q^2\right)(t+q^4t^2)+q^{12}t^3,\cdots\right)$$

be their natural $q-$analogues (cf. [8]).

In previous papers ([4],[5],[6]) we considered the polynomials

$$\left(c_n(t)\right)_{n\geq 1}=\left(1,1+t,1+t+t^2,1+2t+2t^2+t^3,1+2t+4t^2+2t^3+t^4,\cdots\right)$$

which are the limits

(3) $$c_n(t)=C_n\left(t;-1\right)$$

of $C_n(t;q)$ for $q\to -1.$ Their coefficients occur in OEIS, A088855, but apparently, the polynomials $c_n(t)$ themselves seem not to have been studied before.

In [6] we showed that the generating functions

(4) $$C(t,z)=\sum_{n\geq 0}C_n(t)z^n$$

and

(5) $$c(t,z)=\sum_{n\geq 0} c_n(t)z^n$$

are related by

(6) $$c(t,z)c(-t,-z)=C\left(t^2,z^2\right).$$

Here we set $C_0(t)=c_0(t)=1$ by definition.

In the present paper we first recall the fact that both $C_n(t)$ and $c_n(t)$ can be interpreted as certain weights of the set $\Delta_{2n}$ of Dyck paths with length $2n$ and then try to obtain properties of the corresponding weights of the set $\Delta_{2n}^{(m)}$ of the paths from $\Delta_{2n}$ which are contained in the strip $0\leq y\leq m$ for some natural number $m$.

Let me first give some background information.

A path in $\mathbb{Z}^2$ consisting of up-steps $U=(1,1)$ and down-steps $D=(1,-1)$ which starts at some point $(i,j)$ with $i,j\in\mathbb{N}=\{0,1,2,\cdots\}$ and never goes below the $x-$axis will be called *up-down path*. The up-down paths of length $2n$ which start at $(0,0)$ and end at $(2n,0)$ are called *Dyck paths of length* $2n$.

For a given sequence $\tau=\left(t_n\right)_{n\geq 0}$ of real numbers $t_n\neq 0$ with $t_0=1$ we define a *weight* $w=w_\tau$ of the paths in the following way: The weight of each up-step $U$ is $w(U)=1,$ the weight of a down-step $D=D_k$ from height $k+1$ to height $k$ is $w\left(D_k\right)=t_k$. The weight of a path $S_1S_2\cdots S_n$ with steps $S_i$ is the product $w\left(S_1S_2\cdots S_n\right)=\prod_{i=1}^{n} w\left(S_i\right)$ and the weight of a set $\left\{P_j\right\}$ of paths is the sum $\sum_j w\left(P_j\right)$ of their weights.

To get some information about these weights we use results about the polynomials which satisfy

(7) $$p_n(x)=xp_{n-1}(x)-t_{n-2}p_{n-2}(x)$$

with $p_{-1}(x)=0$ and $p_0(x)=1.$

By a theorem of Favard these polynomials are orthogonal with respect to the linear functional $L$ defined by

(8) $$L\left(p_n\right)=[n=0].$$

Since the polynomials $p_n$ have degree $n$ there exist uniquely defined numbers $c(n,k)$ such that

(9) $$x^n=\sum_{k=0}^{n} c(n,k)p_k(x).$$

The numbers $c(n,k)$ satisfy

$$c(n,k)=c(n-1,k-1)+t_k c(n-1,k+1) \tag{10}$$

with $c(0,k)=[k=0]$ and $c(n,k)=0$ for $k<0$.

This follows from

$$\sum_{k=0}^{n} c(n,k)p_k(x)=x^n=x\cdot x^{n-1}=\sum_{k=0}^{n-1} c(n-1,k)xp_k(x)=\sum_{k=0}^{n-1} c(n-1,k)\left(p_{k+1}(x)+t_{k-1}p_{k-1}(x)\right)$$
$$=\sum_{k=0}^{n} p_k(x)\left(c(n-1,k-1)+t_k c(n-1,k+1)\right).$$

Therefore, $c(n,k)$ can be interpreted as the weight of the set of up-down paths from $(0,0)$ to $(n,k)$.

By (9) we get

$$\mu_n=L\left(x^n\right)=c(n,0), \tag{11}$$

i.e. $\mu_n=L\left(x^n\right)$ is the weight of the set of paths of length $n$ which start and end at height $0$.

The generating function $\sum_{n\geq 0}\mu_n z^n$ can be obtained from the following Lemma which is a special case of [1], Theorem 5.8.2.

**Lemma**

*Let $\left(t_n\right)_{n\geq 0}$ with $t_n\neq 0$ and $t_0=1$ be a sequence of real numbers and define $p_n(x)$ by (7) and the linear functional $L$ by (8). Then the generating function $\sum_{n\geq 0}\mu_n z^n$ of the moments $\mu_n=L\left(x^n\right)$ is given by the continued fraction*

$$\sum_{n\geq 0}\mu_n z^n=\cfrac{1}{1-\cfrac{t_0 z^2}{1-\cfrac{t_1 z^2}{1-\ddots}}}=\sum_{n\geq 0}\mu_{2n}z^{2n}. \tag{12}$$

With other words, the generating function of the weights $\mu_{2n}=C_n\left(\tau\right)$ of the set $\Delta_{2n}$ of Dyck paths of length $2n$ is given by

$$C\left(\tau,z\right)=\sum_{n\geq 0}C_n(\tau)z^n=\cfrac{1}{1-\cfrac{t_0 z}{1-\cfrac{t_1 z}{1-\ddots}}}=\frac{1}{1-}\frac{t_0 z}{1-}\frac{t_1 z}{1-}\cdots. \tag{13}$$

We are especially interested in the weights

(14) $$C_n^{(m)}(\tau)=w\left(\Delta_{2n}^{(m)}\right)$$

of the set $\Delta_{2n}^{(m)}$ of the Dyck paths of length $2n$ whose heights are $\leq m$ for some natural number $m$.

An up-down path of length $n$ with height $\leq m$ starting at $(0,0)$ is uniquely determined by the sequence $\left(h_0,h_1,\cdots,h_n\right)$ where $h_i$ is its height at point $i$. Note that $h_0=0, \quad 0\leq h_i\leq m,$ and $\left|h_{i+1}-h_i\right|=1$.

Therefore, it can be interpreted as a walk of length $n$ on the path graph $P_{m+1}$ with vertices $0,1,\cdots,m$. If we define a weight on the edges of $P_{m+1}$ by $w\left(i\to i+1\right)=1$ and $w\left(i+1\to i\right)=t_i$ then $C_n^{(m)}(\tau)$ is the total weight of all walks of length $2n$ on $P_{m+1}$ starting and ending on its initial vertex $0$.

To compute $C_n^{(m)}\left(\tau\right)$ we can consider the matrix

(15) $$A^{(m)}(\tau)=\left(a_{i,j}(\tau)\right)_{i,j=0}^{m}$$

with $a_{i,j}(\tau)=w\left(i\to j\right)$ for $\left|i-j\right|=1$ and $a_{i,j}(\tau)=0$ else. Then the entries $a_{n,i,j}^{(m)}\left(\tau\right)$ of $\left(A^{(m)}(\tau)\right)^n$ are the weight of all walks of length $n$ from $i$ to $j$ in $P_{m+1}$. Thus we get

(16) $$C_n^{(m)}\left(\tau\right)=\left(A^{(m)}(\tau)\right)^{2n}\left(0,0\right).$$

This result depends only on the entries with even subscripts. Therefore we can also consider the following smaller matrices.

For even $m=2k$ we set

(17) $$B^{(2k)}(\tau)=\left(a_{2,2i,2j}^{(2k)}\right)_{i,j=0}^{k}.$$

Then

(18) $$C_n^{(2k)}\left(\tau\right)=\left(B^{(2k)}(\tau)\right)^n(0,0).$$

For odd $m=2k+1$ we define

(19) $$B^{(2k+1)}(\tau)=\left(a_{2,2i,2j}^{(2k+1)}\right)_{i,j=0}^{k}.$$

Then we get

(20) $$C_n^{(2k+1)}\left(\tau\right)=\left(B^{(2k+1)}(\tau)\right)^n(0,0).$$

The generating function $C^{(m)}(\tau,z)=\sum_{n\geq 0}C_n^{(m)}(\tau)z^n$ follows from (13) by setting $t_m=0.$

Thus

(21) $$C^{(m)}(\tau,z)=\sum_{n\geq 0}C_n^{(m)}(\tau)z^n=\frac{1}{1-}\frac{t_0z}{1-}\cdots\frac{t_{m-1}z}{1}.$$

To compute $C^{(m)}(\tau,z)$ we use the following well-known result about continued fractions.

Let

(22) $$\frac{a_1}{b_1+}\frac{a_2}{b_2+}\frac{a_3}{b_3+}\cdots:=\frac{a_1}{b_1+\cfrac{a_2}{b_2+\cfrac{a_3}{b_3+\ddots}}}.$$

The formal convergents

(23) $$C_n=\frac{A_n}{B_n}=\frac{a_1}{b_1+}\frac{a_2}{b_2+}\cdots\frac{a_n}{b_n}$$

are given by

(24) $$\begin{aligned}A_n&=b_nA_{n-1}+a_nA_{n-2},\ A_{-1}=1,\ A_0=0,\\ B_n&=b_nB_{n-1}+a_nB_{n-2},\ \mathrm{B}_{-1}=0,\ B_0=1.\end{aligned}$$

The proof is by induction using

$\frac{A_n}{B_n}=\frac{b_nA_{n-1}+a_nA_{n-2}}{b_nB_{n-1}+a_nB_{n-2}}$ and observing that

$$\frac{A_{n+1}}{B_{n+1}}=\frac{\left(b_n+\frac{a_{n+1}}{b_{n+1}}\right)A_{n-1}+a_nA_{n-2}}{\left(b_n+\frac{a_{n+1}}{b_{n+1}}\right)B_{n-1}+a_nB_{n-2}}=\frac{b_nb_{n+1}A_{n-1}+a_{n+1}A_{n-1}+a_nb_{n+1}A_{n-2}}{b_nb_{n+1}B_{n-1}+a_{n+1}B_{n-1}+a_nb_{n+1}B_{n-2}}$$
$$=\frac{b_{n+1}\left(b_nA_{n-1}+a_nA_{n-2}\right)+a_{n+1}A_{n-1}}{b_{n+1}\left(b_nB_{n-1}+a_nB_{n-2}\right)+a_{n+1}B_{n-1}}=\frac{b_{n+1}A_n+a_{n+1}A_{n-1}}{b_{n+1}B_n+a_{n+1}B_{n-1}}.$$

For the continued fraction (13) we write $A_n=p_n(\tau,z),\ B_n=q_n(\tau,z)$ and get

(25) $$\begin{aligned}p_n(\tau,z)&=p_{n-1}(\tau,z)-t_{n-2}zp_{n-2}(\tau,z)\ \text{ with }\ p_0(\tau,z)=0\text{ and }p_1(\tau,z)=1,\\ q_n(\tau,z)&=q_{n-1}(\tau,z)-t_{n-2}zq_{n-2}(\tau,z)\ \text{ with }\ q_0(\tau,z)=1\text{ and }q_1(\tau,z)=1\end{aligned}$$

and

(26) $$\frac{1}{1-}\frac{t_0z}{1-}\cdots\frac{t_{m-2}z}{1-}\frac{t_{m-1}z}{1}=\frac{p_{m+1}(\tau,z)}{q_{m+1}(\tau,z)}=C^{(m)}(\tau,z),$$

We set $C^{(0)}(\tau,z)=1$ since the only Dyck path of height $0$ is the trivial path of length $0$ which by definition has weight $1$.

For the constant sequence $\tau=(1,1,1,\cdots)$ the weight of paths reduces to their number and (13) reduces to the well-known fact that the generating function

(27) $$C(z)=\sum_{n\geq 0} C_n z^n = \frac{1-\sqrt{1-4z}}{2z}$$

can be represented by the corresponding continued fraction

(28) $$C(z)=\cfrac{1}{1-\cfrac{z}{1-\cfrac{z}{1-\ddots}}}=\frac{1}{1-}\frac{z}{1-}\frac{z}{1-}\cdots.$$

**2. Weights of bounded Dyck paths which are related to Narayana polynomials**

As shown in [6] the Narayana polynomials

(29) $$C_n(t)=\sum_{k=0}^{n-1}\binom{n}{k}\binom{n-1}{k}\frac{t^k}{k+1}$$

can be interpreted as the weight

(30) $$C_n(t)=w_+(\Delta_{2n})$$

if $w_+$ is the weight corresponding to the sequence

(31) $$\tau_+(t)=(1,t,1,t,1,t,\cdots).$$

Thus the generating function

(32) $$C(t,z)=\sum_{n\geq 0} C_n(t) z^n$$

has the continued fraction expansion

(33) $$C(t,z)=\frac{1}{1-}\frac{z}{1-}\frac{tz}{1-}\frac{z}{1-}\frac{tz}{1-}\cdots=\cfrac{1}{1-\cfrac{z}{1-tzC(t,z)}}.$$

Analogously the polynomials

(34) $$c_n(t)=C_n(t;-1)$$

can be interpreted as the weight

(35) $$c_n(t)=w_-(\Delta_{2n})$$

if $w_-$ corresponds to the sequence

(36) $$\tau_-(t) = (1, t, -1, -t, 1, t, -1, -t, \cdots).$$

As shown in [4]

(37) $$c_n(t) = \sum_{k=0}^{n} \left\langle {n \atop k} \right\rangle t^k$$

with

(38) $$\left\langle {n \atop k} \right\rangle = \binom{\left\lfloor \frac{n-1}{2} \right\rfloor}{\left\lfloor \frac{k}{2} \right\rfloor} \binom{\left\lfloor \frac{n}{2} \right\rfloor}{\left\lfloor \frac{k+1}{2} \right\rfloor}.$$

The generating function is given by

(39) $$c(t,z) = \sum_{n \geq 0} c_n(t) z^n = \frac{1}{1-} \frac{z}{1-} \frac{tz}{1+} \frac{z}{1+} \frac{tz}{1-} \cdots.$$

We are especially interested in concrete facts about the polynomials

(40) $$C_n^{(m)}(t) = C_n^{(m)}(\tau_+(t))$$

and

(41) $$c_n^{(m)}(t) = C_n^{(m)}(\tau_-(t)).$$

By (15) we get

(42) $$A_+^{(m)}(t) := A^{(m)}(\tau_+(t)) = \left(a_{i,j}(t)\right)_{i,j=0}^{m}$$

with $a_{j+1,j}(t) = 1$ for even $j$, $a_{j+1,j}(t) = t$ for odd $t$, $a_{i,i+1}(t) = 1$ and $a_{i,j}(t) = 0$ else.

By (17) we get

(43) $$B_+^{(2k)}(t) := B^{(2k)}(\tau_+(t)) = \left(u_{2i,2j}^{(2k)}\right)_{i,j=0}^{k}$$

where $u_{2i,2j}^{(2k)}$ is the weight of the up-down paths of length 2 from a point with height $2i$ to a point with height $2j$. Thus $u_{2i+2,2i}^{(2k)} = w_+(D_{2i+1}D_{2i}) = t$, $u_{2i,2i+2}^{(2k)} = w(UU) = 1$, $u_{0,0}^{(2k)} = w_+(UD_0) = 1$, $u_{2k,2k}^{(2k)} = w_+(D_{2k-1}U) = t$, and for $0 < i < k$ $u_{2i,2i}^{(2k)} = w_+(UD_{2i}) + w_+(D_{2i-1}U) = 1 + t$ and $u_{2i,2j}^{(2k)} = 0$ else.

For example, $B_+^{(2)}(t) = \begin{pmatrix} 1 & 1 \\ t & t \end{pmatrix}$, $B_+^{(4)}(t) = \begin{pmatrix} 1 & 1 & 0 \\ t & 1+t & 1 \\ 0 & t & t \end{pmatrix}$.

By (19) we get

(44) $$B_+^{(2k+1)}(t) := B^{(2k+1)}\left(\tau_+(t)\right) = \left(u_{2i,2j}^{(2k+1)}\right)_{i,j=0}^{k}$$

with $u_{2i,2j}^{(2k+1)} = u_{2i,2j}^{(2k)}$ with the exception that $u_{2k,2k}^{(2k+1)} = w_+\left(UD_{2k}\right) + w_+\left(D_{2k-1}U\right) = 1+t.$

For example,

$$B_+^{(3)}(t) = \begin{pmatrix} 1 & 1 \\ t & 1+t \end{pmatrix},\ B_+^{(5)}(t) = \begin{pmatrix} 1 & 1 & 0 \\ t & 1+t & 1 \\ 0 & t & 1+t \end{pmatrix}.$$

**Theorem 1**

*The polynomials* $C_n^{(m)}(t)$ *are given by*

(45) $$C_n^{(2k)}(t) = \left(B_+^{(2k)}\right)^n (0,0), \quad C_n^{(2k+1)}(t) = \left(B_+^{(2k+1)}\right)^n (0,0).$$

For example,

$$\left(\left(B_+^{(2)}(t)\right)^n\right)_{n\ge1} = \left(\begin{pmatrix} 1 & 1 \\ t & t \end{pmatrix}, \begin{pmatrix} 1+t & 1+t \\ t+t^2 & t+t^2 \end{pmatrix}, \begin{pmatrix} 1+2t+t^2 & 1+2t+t^2 \\ t+2t^2+t^3 & t+2t^2+t^3 \end{pmatrix}, \begin{pmatrix} 1+3t+3t^2+t^3 & 1+3t+3t^2+t^3 \\ t+3t^2+3t^3+t^4 & t+3t^2+3t^3+t^4 \end{pmatrix}, \cdots\right)$$

and $C_{n+1}^{(2)}(t) = (1+t)^n$.

$$\left(\left(B_+^{(3)}(t)\right)^n\right)_{n\ge1} = \left(\begin{pmatrix} 1 & 1 \\ t & 1+t \end{pmatrix}, \begin{pmatrix} 1+t & 2+t \\ 2t+t^2 & 1+3t+t^2 \end{pmatrix}, \begin{pmatrix} 1+3t+t^2 & 3+4t+t^2 \\ 3t+4t^2+t^3 & 1+6t+5t^2+t^3 \end{pmatrix}, \begin{pmatrix} 1+6t+5t^2+t^3 & 4+10t+6t^2+t^3 \\ 4t+10t^2+6t^3+t^4 & 1+10t+15t^2+7t^3+t^4 \end{pmatrix}\right)$$

and $\left(C_{n+1}^{(3)}(t)\right)_{n\ge0} = \left(1, 1+t, 1+3t+t^2, 1+6t+5t^2+t^3, \cdots\right).$

Let us now state the analogous facts for $c_n^{(m)}(t)$. Here we have

(46) $$A_-^{(m)}(t) := A^{(m)}\left(\tau_-(t)\right) = \left(a_{i,j}(t)\right)_{i,j=0}^{m}$$

with $a_{2j+1,2j}(t) = (-1)^j$, $a_{2j+2,2j+1}(t) = (-1)^j t$, $a_{i,i+1}(t) = 1$ and $a_{i,j}(t) = 0$ else.

(47) $$B_-^{(2k)}(t) = B^{(2k)}\left(\tau_-(t)\right) = \left(b_{2i,2j}^{(2k)}(t)\right)_{i,j=0}^{k}$$

with $b_{0,0}^{(2k)}(t) = 1$, $b_{2k,2k}^{(2k)}(t) = (-1)^{k-1}t$, $b_{2i,2i}^{(2k)}(t) = (-1)^i(1-t)$ for $1 \le i \le k-1$, $b_{2i,2i-2}^{(2k)} = t$, $b_{2i,2i+2}^{(2k)} = 1$, $b_{2i,2j}^{(2k)} = 0$ else.

For example,

$$B_-^{(4)}(t)=\begin{pmatrix}1&1&0\\t&-1+t&1\\0&t&-t\end{pmatrix},\quad B_-^{(6)}(t)=\begin{pmatrix}1&1&0&0\\t&-1+t&1&0\\0&t&1-t&1\\0&0&t&t\end{pmatrix}.$$

(48)
$$B_-^{(2k+1)}(t):=B^{(2k+1)}\left(\tau_-(t)\right)=\left(u_{2i,2j}^{(2k+1)}\right)_{i,j=0}^{k}$$

with $u_{2i,2j}^{(2k)}=u_{2i,2j}^{(2k+1)}$ except that $u_{2k,2k}^{(2k+1)}=(t-1)(-1)^{k-1}$.

For example,

$$B_-^{(5)}(t)=\begin{pmatrix}1&1&0\\t&-1+t&1\\0&t&1-t\end{pmatrix},\ B_-^{(7)}(t)=\begin{pmatrix}1&1&0&0\\t&-1+t&1&0\\0&t&1-t&1\\0&0&t&-1+t\end{pmatrix}.$$

**Theorem 2**

*The polynomials $c_n^{(m)}(t)$ are given by*

(49)
$$c_n^{(2k)}(t)=\left(B_-^{(2k)}\right)^n(0,0),$$
$$c_n^{(2k+1)}(t)=\left(B_-^{(2k+1)}\right)^n(0,0).$$

For example,

$$\left(\left(B_-^{(3)}\right)^n\right)_{n\ge1}=\left(\begin{pmatrix}1&1\\t&-1+t\end{pmatrix},\begin{pmatrix}1+t&t\\t^2&1-t+t^2\end{pmatrix},\begin{pmatrix}1+t+t^2&1+t^2\\t+t^3&-1+2t-t^2+t^3\end{pmatrix},\begin{pmatrix}1+2t+t^2+t^3&2t+t^3\\2t^2+t^4&1-2t+3t^2-t^3+t^4\end{pmatrix},\cdots\right)$$

and $\left(c_n^{(3)}(t)\right)_{n\ge1}=\left(1,1+t,1+t+t^2,1+2t+t^2+t^3,1+2t+3t^2+t^3+t^4,\cdots\right)$.

We now consider the generating functions

(50)
$$C^{(m)}\left(\tau,z\right)=\sum_{n\ge0}C_n^{(m)}(\tau)z^n=\frac{1}{1-}\frac{t_0z}{1-}\cdots\frac{t_{m-1}z}{1}.$$

for $\tau=\tau_+(t)$ and $\tau=\tau_-(t)$.

We shall use the notations

(51)
$$C^{(m)}\left(\tau_+(t),z\right)=\frac{P^{(m+1)}(t,z)}{Q^{(m+1)}(t,z)}=C^{(m)}(t,z)=\sum_{n\ge0}C_n^{(m)}(t)z^n$$

and

(52) $$C^{(m)}\left(\tau_-(t),z\right)=\frac{p^{(m+1)}(t,z)}{q^{(m+1)}(t,z)}=c^{(m)}(t,z)=\sum_{n\geq 0}c_n^{(m)}(t)z^n.$$

The first terms are

$$\left(C^{(m)}(t,z)\right)_{m\geq 0}=\left(1,\frac{1}{1-z},\frac{1-tz}{1-(1+t)z},\frac{1-(1+t)z}{1-(2+t)z+z^2},\frac{1-(1+2t)z+t^2z^2}{1-2(1+t)z+\left(1+t+t^2\right)z^2},\cdots\right)$$

$$\left(c^{(m)}(t,z)\right)_{m\geq 0}=\left(1,\frac{1}{1-z},\frac{1-tz}{1-(1+t)z},\frac{1+(1-t)z}{1-tz-z^2},\frac{1+z-t^2z^2}{1-\left(1+t+t^2\right)z^2},\cdots\right)$$

We get

$$\deg Q^{(m+1)}(t,z)=\deg q^{(m+1)}(t,z)=\left\lfloor\frac{m+1}{2}\right\rfloor,$$
$$\deg P^{(m+1)}(t,z)=\deg p^{(m+1)}(t,z)=\left\lfloor\frac{m}{2}\right\rfloor.$$

This follows from

$$\frac{P^{(m+1)}(t,z)}{Q^{(m+1)}(t,z)}=C^{(m)}(t,z)=\frac{1}{1-\dfrac{z}{1-tzC^{(m-2)}(t,z)}}=\frac{1-tzC^{(m-2)}(t,z)}{1-z-tzC^{(m-2)}(t,z)}$$

$$=\frac{1-tz\dfrac{P^{(m-1)}(t,z)}{Q^{(m-1)}(t,z)}}{1-z-tz\dfrac{P^{(m-1)}(t,z)}{Q^{(m-1)}(t,z)}}=\frac{Q^{(m-1)}(t,z)-tzP^{(m-1)}(t,z)}{Q^{(m-1)}(t,z)(1-z)+tzP^{(m-1)}(t,z)}$$

by induction.

It should be noted that $C_n^{(m)}(t)=C_n(t)$ and $c_n^{(m)}(t)=c_n(t)$ for $n\leq m.$

**3. The case t=1**

For $t_i=1$ let $C_n^{(m)}$ be the number of Dyck paths of length $2n$ with height $\leq m.$

Here $A_+^{(m)}:=A^{(m)}\left(\tau_+(1)\right)=\left(a_{i,j}\right)_{i,j=0}^m$ with $a_{i,j}=1$ for $|i-j|=1$ and $a_{i,j}=0$ else and

$B_+^{(2k)}:=B^{(2k)}\left(\tau_+(1)\right)=\left(b_{i,j}\right)_{i,j=0}^k$ with $b_{0,0}=b_{k,k}=1,$ $b_{i,i}=2$ for $1\leq i\leq k-1,$ $b_{i,j}=1$ for $|i-j|=1$ and $b_{i,j}=0$ else.

$B_+^{(2k+1)}=\left(\beta_{i,j}\right)_{i,j=0}^k$ with $\beta_{i,j}=b_{i,j}$ except $\beta_{2k,2k}=2.$

To obtain closed formulas for generating functions we use some well-known results about the Fibonacci polynomials

(53) $$F_n(x,s)=\sum_{j=0}^{\left\lfloor\frac{n-1}{2}\right\rfloor}\binom{n-1-j}{j}s^j x^{n-1-2j}$$

which satisfy

(54) $$F_n(x,s)=xF_{n-1}(x,s)+sF_{n-2}(x,s)$$

with initial values $F_0(x,s)=0$ and $F_1(x,s)=1$ and the Lucas polynomials which satisfy

(55) $$L_n(x,s)=xL_{n-1}(x,s)+sL_{n-2}(x,s)$$

with initial values $L_0(x,s)=2$ and $L_1(x,s)=x.$ Binet's formulae give

(56) $$F_n(x,s)=\frac{\alpha(x,s)^n-\beta(x,s)^n}{\alpha(x,s)-\beta(x,s)}$$

and

(57) $$L_n(x,s)=\alpha(x,s)^n+\beta(x,s)^n=F_{n+1}(x,s)+sF_{n-1}(x,s)$$

with

(58) $$\alpha(x,s)=\frac{x+\sqrt{x^2+4s}}{2},\quad \beta(x,s)=\frac{x-\sqrt{x^2+4s}}{2}.$$

Note that $\alpha(x,s)^2=x\alpha(x,s)+s,\ \beta(x,s)^2=x\beta(x,s)+s,$ and $\alpha(x,s)\beta(x,s)=-s.$

We set $P^{(m)}(1,z)=P^{(m)}(z),\ Q^{(m)}(1,z)=Q^{(m)}(z),\ p^{(m)}(1,z)=p(z),\ q^{(m)}(1,z)=q^{(m)}(z).$

It is well known that

(59) $$C^{(m)}(z)=\frac{P_{m+1}(z)}{Q_{m+1}(z)}=\frac{F_{m+1}(1,-z)}{F_{m+2}(1,-z)}.$$

This also follows from (54) by comparing with (25) which reduces to

$P_n(z)=P_{n-1}(z)-zP_{n-2}(z)$ with $P_0(z)=0$ and $P_1(z)=1,$
$Q_n(z)=Q_{n-1}(z)-zQ_{n-2}(z)$ with $Q_0(z)=1$ and $Q_1(z)=1.$

Let us consider the sequences $\left(C_n^{(m)}\right)_{n\geq 0}$ for small $m.$

$\left(C_n^{(0)}\right)=(1,0,0,0,\cdots),$

$\left(C_n^{(1)}\right)=(1,1,1,\cdots)$,

$\left(C_n^{(2)}\right)=(1,1,2,4,8,16,\cdots)$, OEIS, A011782,

$\left(C_n^{(3)}\right)=(1,1,2,5,13,34,89,233,\cdots)$, OEIS, A 001519, i.e., $C_n^{(3)}=F_{2n-1}$ where

$(F_n)_{n\geq 0}=(0,1,1,2,3,5,8,\cdots)$ is the sequence of Fibonacci numbers and $F_{-1}=1$.

$\left(C_n^{(4)}\right)=(1,1,2,5,14,41,122,365,1094,\cdots)$, OEIS, A124302, i.e., $C_{n+1}^{(4)}=\frac{3^n+1}{2}$.

For general $m$ we have $C_n^{(m)}=C_n$ for $n\leq m$ and $F_{m+2}(1,-z)\sum_{n\geq 0}C_n^{(m)}z^n=F_{m+1}(1,-z)$.

Let us also consider $C^{(m)}(-1,z)=\frac{p_{m+1}(-1,z)}{q_{m+1}(-1,z)}$.

By (25) we get $p_m(-1,z)=p_{m-1}(-1,z)-(-1)^m z p_{m-2}(-1,z)$ with $p_0(-1,z)=0,\ \ p_1(-1,z)=1$.

This gives

$(p_m(-1,z))_{m\geq 0}=\left(0,1,1,1+z,1,1+z+z^2,1+z^2,1+z+2z^2+z^3,1+2z^2,1+z+3z^2+2z^3+z^4,1+3z^2+z^4,\cdots\right)$.

Comparing with

$\left(F_m\left(1,z^2\right)\right)=\left(1,1,1+z^2,1+2z^2,1+3z^2+z^4,1+4z^2+3z^4,1+5z^2+6z^4+z^6,\cdots\right)$

we see that

(60)
$$\begin{aligned}p_{2m}(-1,z)&=F_m\left(1,z^2\right),\\ p_{2m-1}(-1,z)&=F_m\left(1,z^2\right)+zF_{m-1}\left(1,z^2\right).\end{aligned}$$

This gives $C^{(2m)}(-1,z)=\frac{F_{m+1}\left(1,z^2\right)+zF_m\left(1,z^2\right)}{F_{m+1}\left(1,z^2\right)}=1+z\frac{F_m\left(1,z^2\right)}{F_{m+1}\left(1,z^2\right)}$ and therefore

(61)
$$C^{(2m)}(-1,z)=1+zC^{(m-1)}\left(-z^2\right).$$

For $C^{(2m-1)}(-1,z)$ we get analogously

(62)
$$C^{(2m-1)}(-1,z)=\frac{F_m\left(1,z^2\right)}{F_{m+1}\left(1,z^2\right)-zF_m\left(1,z^2\right)}.$$

More interesting are the generating functions $c^{(m)}(z)$.

$$\left(c^{(m)}(z)\right)_{m\geq 0}=\left(1,\frac{1}{1-z},\frac{1-z}{1-2z},\frac{1}{1-z-z^2},\frac{1+z-z^2}{1-3z^2},\frac{(1-z)(1+z)}{1-z-2z^2+z^3},\frac{1-z-2z^2+z^3}{(1-2z)\left(1-2z^2\right)},\cdots\right)$$

The corresponding sequences $\left(c_n^{(m)}\right)_{n\geq 0}$ are

$\left(c_n^{(1)}\right)=(1,1,1,\cdots)$,

$\left(c_n^{(2)}\right)=(1,1,2,4,8,16,\cdots)$,

$\left(c_n^{(3)}\right)=(1,1,2,3,5,8,13,21,34,\cdots)$, i.e., $c_n^{(3)}=F_{n+1}$,

$\left(c_n^{(4)}\right)=(1,1,2,3,6,9,18,27,54,\cdots)$, OEIS, A182522, i.e., $c_{2n+1}^{(3)}=3^n$ and $c_{2n+2}^{(3)}=2\cdot 3^n$.

$\left(c_n^{(5)}\right)=(1,1,2,3,6,10,19,33,61,108,197,352,638,\cdots)$, OEIS, A028495.

These sequences satisfy $c_n^{(m)}=c_n=B_n=\binom{n}{\left\lfloor \frac{n}{2}\right\rfloor}$ for $n\leq m$.

Let us write $c^{(m)}(z)=\dfrac{p_{m+1}(z)}{q_{m+1}(z)}$.

By (25) we get

$p_n(z)=p_{n-1}(z)-(-1)^{\binom{n-2}{2}}zp_{n-2}(z)$ with $p_0(z)=0$ and $p_1(z)=1$,

$q_n(z)=q_{n-1}(z)-(-1)^{\binom{n-2}{2}}zq_{n-2}(z)$ with $q_0(z)=1$ and $q_1(z)=1$.

This gives

$$\begin{aligned}
p_{2n}(z) &= p_{2n-1}(z)-(-1)^{\binom{2n-2}{2}}zp_{2n-2}(z)=p_{2n-1}(z)+(-1)^n zp_{2n-2}(z)\\
&= p_{2n-2}(z)-(-1)^{\binom{2n-3}{2}}zp_{2n-3}(z)+(-1)^n zp_{2n-2}(z)\\
&= p_{2n-2}(z)\left(1+(-1)^n z\right)-(-1)^n zp_{2n-3}(z)\\
&= p_{2n-2}(z)\left(1+(-1)^n z\right)+(-1)^{n-1}z\left(p_{2n-2}(z)+(-1)^n zp_{2n-4}\right)\\
&= p_{2n-2}(z)-z^2p_{2n-4}(z)
\end{aligned}$$

Since $p_0(z)=0$ and $p_1(z)=1$ we get $p_2(z)=1$ and thus

(63) $$p_{2n}(z)=F_n\left(1,-z^2\right).$$

$q_0(z)=q_1(z)=1$ gives $q_2(z)=1-z$ and therefore

(64) $$q_{2n}(z)=F_{n+1}\left(1,-z^2\right)-zF_n\left(1,-z^2\right).$$

From $p_{2n}(z)=p_{2n-1}(z)+(-1)^n zp_{2n-2}(z)$ and (63) we get

(65) $$p_{2n+1}(z)=F_{n+1}\left(1,-z^2\right)+(-1)^n zF_n\left(1,-z^2\right).$$

Analogously we get $q_{2n+1}(z)=q_{2n+2}(z)+(-1)^n zq_{2n}(z)$. This gives

$$q_{4n+1}(z)=q_{4n+2}(z)+zq_{4n}(z)=F_{2n+2}\left(1,-z^2\right)-zF_{2n+1}\left(1,-z^2\right)+z\left(F_{2n+1}\left(1,-z^2\right)-zF_{2n}\left(1,-z^2\right)\right)$$
$$=F_{2n+2}\left(1,-z^2\right)-z^2F_{2n}\left(1,-z^2\right).$$

By (57) this gives

(66) $$q_{4n+1}(z)=L_{2n+1}\left(1,-z^2\right).$$

$$q_{4n+3}(z)=q_{4n+4}(z)-zq_{4n+2}(z)=F_{2n+3}\left(1,-z^2\right)-zF_{2n+2}\left(1,-z^2\right)$$
$$-z\left(F_{2n+2}\left(1,-z^2\right)-zF_{2n+1}\left(1,-z^2\right)\right)=F_{2n+3}\left(1,-z^2\right)-2zF_{2n+2}\left(1,-z^2\right)+z^2F_{2n+1}\left(1,-z^2\right)$$
$$=(1-2z)F_{2n+2}\left(1,-z^2\right)$$

gives

(67) $$q_{4n+3}(z)=(1-2z)F_{2n+2}\left(1,-z^2\right).$$

**Theorem 3**

(68) $$c^{(2m-1)}(z)=\frac{p_{2m}(z)}{q_{2m}(z)}=\frac{F_m\left(1,-z^2\right)}{F_{m+1}\left(1,-z^2\right)-zF_m\left(1,-z^2\right)},$$

(69) $$c^{(4m)}(z)=\frac{p_{4m+1}(z)}{q_{4m+1}(z)}=\frac{F_{2m+1}\left(1,-z^2\right)+zF_{2m}\left(1,-z^2\right)}{L_{2m+1}\left(1,-z^2\right)},$$

(70) $$c^{(4m+2)}(z)=\frac{p_{4m+3}(z)}{q_{4m+3}(z)}=\frac{F_{2m+2}\left(1,-z^2\right)-zF_{2m+1}\left(1,-z^2\right)}{(1-2z)F_{2m+2}\left(1,-z^2\right)}.$$

**Corollary**

*For* $m\not\equiv 2 \bmod 4$

(71) $$c^{(m)}(z)=\frac{r_m(z)}{F_{m+2}\left(1,-z^2\right)}$$

*for some polynomial* $r_m(z)$.

**Proof**

$$q_{2m}(z)q_{2m}(-z)=\left(F_{m+1}\left(1,-z^2\right)-zF_m\left(1,-z^2\right)\right)\left(F_{m+1}\left(1,-z^2\right)+zF_m\left(1,-z^2\right)\right)$$
$$=F_{m+1}\left(1,-z^2\right)^2-z^2F_m\left(1,-z^2\right)^2=F_{2m+1}\left(1,-z^2\right).$$

With $\alpha = \alpha\left(1, -z^2\right)$ and $\beta = \beta\left(1, -z^2\right)$ the last identity is

$$\left(\frac{\alpha^{m+1} - \beta^{m+1}}{\alpha - \beta}\right)^2 - z^2\left(\frac{\alpha^m - \beta^m}{\alpha - \beta}\right)^2 = \left(\frac{\alpha^{2m+1} - \beta^{2m+1}}{\alpha - \beta}\right)$$

which is easily verified using $\alpha\beta = z^2$. Thus,

$$c^{(2m-1)}(z) = \frac{p_{2m}(z)}{q_{2m}(z)} = \frac{p_{2m}(z)q_{2m}(-z)}{q_{2m}(z)q_{2m}(-z)} = \frac{p_{2m}(z)q_{2m}(-z)}{F_{2m+1}\left(1, -z^2\right)}.$$

For $c^{(4m)}(z)$ we get

$$q_{4m+1}(z)F_{2m+1}(1, -z^2) = L_{2m+1}(1, -z^2)F_{2m+1}(1, -z^2) = \left(\alpha^{2m+1} + \beta^{2m+1}\right)\frac{\left(\alpha^{2m+1} - \beta^{2m+1}\right)}{\alpha - \beta}$$

$$= \frac{\alpha^{4m+2} - \beta^{4m+2}}{\alpha - \beta} = F_{2m+2}(1, -z^2)$$

which implies $c^{(4m)}(z) = \dfrac{p_{4m+1}(z)F_{2m+1}(1, -z^2)}{q_{4m+1}(z)F_{2m+1}(1, -z^2)} = \dfrac{p_{4m+1}(z)F_{2m+1}(1, -z^2)}{F_{4m+2}\left(1, -z^2\right)}.$

In [2] we studied the numbers $\alpha(n, k)$

(72) $$\alpha(n, k) = \sum_{j \in \mathbb{Z}} (-1)^j \binom{n}{\left\lfloor \frac{n + (k+2)j}{2} \right\rfloor},$$

which can be interpreted as the number of all up-down paths from $(0, 0)$ to $(n, 0)$ which are contained in the strip $0 \le y \le k$.

By comparing the generating functions $c^{(m)}(z)$ with formulas [2], (1.25) and (1.26) we get

**Theorem 4**

*For $m \not\equiv 2 \bmod 4$ we get*

*(73)* $$c_n^{(m)} = \alpha(n, m) = \sum_{j=0}^{m} \left(A_+^{(m)}\right)^n (0, j).$$

*i.e. $c_n^{(m)}$ is the number of up-down paths of length $n$ which are $\le m$.*

Some remarks in OEIS, A007582, A085282 and computations for small $m \equiv 2 \bmod 4$ suggest

**Conjecture 1**

*Let* $m = 4k-2$ *and* $U_k = \left(u_{i,j}\right)_{i,j=0}^{4k-1}$ *with* $u_{i,j} = 1$ *for* $j - i \equiv \pm 1 \bmod (4k)$ *and* $u_{i,j} = 0$ *else. Then*

(74)
$$c_{2n+1}^{(4k-2)} = U_k^{2n+1}(0,1),$$
$$c_{2n+2}^{(4k-2)} = U_k^{2n+2}(0,0) = 2c_{2n+1}^{(4k-2)}.$$

## 4. Analogues of Pascal's triangle

For $m \geq 2$ the polynomials $C_{n+1}^{(m)}(t)$ and $c_{n+1}^{(m)}(t)$ have degree $n$ because the path $U\left(UD_1\right)^n D_0$ is the only path in $\Delta_{2(n+1)}$ with weight $t^n$ and the weight of other paths is a multiple of $t^i$ with $i < n$. For $m = 2$ we have $C^{(2)}(t,z) = c^{(2)}(t,z) = \dfrac{1-tz}{1-(1+t)z} = 1 + \dfrac{z}{1-(1+t)z}$ with $C_{n+1}^{(2)}(t) = c_{n+1}^{(2)}(t) = (1+t)^n = \sum_{k=0}^{n} \binom{n}{k} t^k$. This suggests to write

(75)
$$C_{n+1}^{(m)}(t) = \sum_{k=0}^{n} \begin{bmatrix} n \\ k \end{bmatrix}^{(m)} t^k,$$
$$c_{n+1}^{(m)}(t) = \sum_{k=0}^{n} \binom{n}{k}^{(m)} t^k$$

and consider the triangles $\left(\begin{bmatrix} n \\ k \end{bmatrix}^{(m)}\right)$ and $\left(\binom{n}{k}^{(m)}\right)$ as analogues of Pascal's triangle.

For small $m$ the corresponding number triangles are known in other contexts. I shall give a short survey with their identifiers in OEIS.

For $m = 3$ we get

$\left(C_{n+1}^{(3)}(t)\right)_{n \geq 0} = \left(1, 1+t, 1+3t+t^2, 1+6t+5t^2+t^3, 1+10t+15t^2+7t^3+t^4, \cdots\right)$, A085478,

with $\begin{bmatrix} n \\ k \end{bmatrix}^{(3)} = \binom{n+k}{2k}$, and

$\left(c_{n+1}^{(3)}(t)\right)_{n \geq 0} = \left(1, 1+t, 1+t+t^2, 1+2t+t^2+t^3, 1+2t+3t^2+t^3+t^4, \cdots\right)$, A046854,

with $\binom{n}{k}^{(3)} = \binom{\left\lfloor \frac{n+k}{2} \right\rfloor}{k} = \begin{bmatrix} n+k \\ 2k \end{bmatrix}_{q=-1}$.

The first terms of $\binom{n}{k}^{(3)}$ are

$$\left(\binom{n}{k}^{(3)}\right)=\begin{pmatrix}1&0&0&0&0&0&0\\1&1&0&0&0&0&0\\1&1&1&0&0&0&0\\1&2&1&1&0&0&0\\1&2&3&1&1&0&0\\1&3&3&4&1&1&0\\1&3&6&4&5&1&1\end{pmatrix}.$$

Since $c^{(3)}(t,z)=\dfrac{1+(1-t)z}{1-tz-z^2}$ we get $c_{n+1}^{(3)}(t)=tc_n^{(3)}(t)+c_{n-1}^{(3)}(t)$ which implies

(76) $$\binom{n}{k}^{(3)}=\binom{n-1}{k-1}^{(3)}+\binom{n-2}{k}^{(3)}.$$

The first terms of $\begin{bmatrix}n\\k\end{bmatrix}^{(3)}$ are

$$\left(\begin{bmatrix}n\\k\end{bmatrix}^{(3)}\right)=\begin{pmatrix}1&0&0&0&0&0&0\\1&1&0&0&0&0&0\\1&3&1&0&0&0&0\\1&6&5&1&0&0&0\\1&10&15&7&1&0&0\\1&15&35&28&9&1&0\\1&21&70&84&45&11&1\end{pmatrix}$$

Since $C^{(3)}(t,z)=\dfrac{1-(1+t)z}{1-2z-tz+z^2}$ we get $C_{n+1}^{(3)}(t)=2C_n^{(3)}(t)+tC_n^{(3)}(t)-C_{n-1}^{(3)}(t)$ which gives

$$\begin{bmatrix}n\\k\end{bmatrix}^{(3)}=2\begin{bmatrix}n-1\\k\end{bmatrix}^{(3)}+\begin{bmatrix}n-1\\k-1\end{bmatrix}^{(3)}-\begin{bmatrix}n-2\\k\end{bmatrix}^{(3)}.$$

For $m=4$ we get

$\left(C_{n+1}^{(4)}(t)\right)_{n\geq 0}=\left(1,1+t,1+3t+t^2,1+6t+6t^2+t^3,1+10t+19t^2+10t^3+t^4,\cdots\right)$, A056241,

with $\sum_{k=0}^{n}\begin{bmatrix}n\\k\end{bmatrix}^{(4)}t^{2k}=\dfrac{\left(1+t+t^2\right)^n+\left(1-t+t^2\right)^n}{2}$, and

$\left(c_{n+1}^{(4)}(t)\right)_{n\geq 0}=\left(1,1+t,1+t+t^2,1+2t+2t^2+t^3,1+2t+3t^2+2t^3+t^4,\cdots\right)$, A169623,

with $\sum_{k=0}^{2n}\binom{2n}{k}^{(4)}t^k=\left(1+t+t^2\right)^n$, $\sum_{k=0}^{2n+1}\binom{2n+1}{k}^{(4)}t^k=(1+t)\left(1+t+t^2\right)^n$.

It should be noted that

$$C^{(4)}\left(t^2,z^2\right)=\frac{1-(1+2t^2)z^2+t^4z^4}{1-2(1+t^2)z^2+\left(1+t^2+t^4\right)z^4}=\frac{\left(1-t^2z^2+z\right)\left(1-t^2z^2-z\right)}{\left(1-z^2-t^2z^2-tz^2\right)\left(1-z^2-t^2z^2+tz^2\right)}$$
$$=c^{(4)}(t,z)c^{(4)}(-t,-z).$$

The first terms of $\binom{n}{k}^{(4)}$ are

$$\left(\binom{n}{k}^{(4)}\right)=\begin{pmatrix}1&0&0&0&0&0&0\\1&1&0&0&0&0&0\\1&1&1&0&0&0&0\\1&2&2&1&0&0&0\\1&2&3&2&1&0&0\\1&3&5&5&3&1&0\\1&3&6&7&6&3&1\end{pmatrix}$$

From $c^{(4)}(t,z)=\dfrac{1+z-t^2z^2}{1-\left(1+t+t^2\right)z^2}$ we get $c_{n+1}^{(4)}(t)=(1+t+t^2)c_{n-1}^{(4)}(t)$ which implies

(77) $$\binom{n}{k}^{(4)}=\binom{n-2}{k}^{(4)}+\binom{n-2}{k-1}^{(4)}+\binom{n-2}{k-2}^{(4)}.$$

$$\left(\begin{bmatrix}n\\k\end{bmatrix}^{(4)}\right)=\begin{pmatrix}1&0&0&0&0&0&0\\1&1&0&0&0&0&0\\1&3&1&0&0&0&0\\1&6&6&1&0&0&0\\1&10&19&10&1&0&0\\1&15&45&45&15&1&0\\1&21&90&141&90&21&1\end{pmatrix}$$

From $C^{(4)}(t,z)=\dfrac{1-z-2tz+t^2z^2}{1-2z-2tz+z^2+tz^2+t^2z^2}$ we get

(78) $$\begin{bmatrix} n \\ k \end{bmatrix}^{(4)} = 2\begin{bmatrix} n-1 \\ k \end{bmatrix}^{(4)} + 2\begin{bmatrix} n-1 \\ k-1 \end{bmatrix}^{(4)} - \begin{bmatrix} n-2 \\ k \end{bmatrix}^{(4)} - \begin{bmatrix} n-2 \\ k-1 \end{bmatrix}^{(4)} - \begin{bmatrix} n-2 \\ k-2 \end{bmatrix}^{(4)}.$$

Both sequences $C_n^{(4)}(t)$ and $c_n^{(4)}(t)$ are unimodal and palindromic. Therefore (cf. [8], 4.4) they have a representation as linear combination of terms $t^j(1+t)^{n-2j}$. Computations suggest that

$$C_{n+1}^{(4)}(t) = \sum_k \begin{bmatrix} n \\ k \end{bmatrix}^{(4)} t^k = \sum_j \binom{n}{2j} t^j (1+t)^{n-2j},$$

$$c_{n+1}^{(4)}(t) = \sum_k \binom{n}{k}^{(4)} t^k = \sum_j (-1)^j \binom{\left\lfloor \frac{n}{2} \right\rfloor}{j} t^j (1+t)^{n-2j}.$$

$$\left(C_{n+1}^{(5)}(t)\right) = \left(1, 1+t, 1+3t+t^2, 1+6t+6t^2+t^3, 1+10t+20t^2+10t^3+t^4, 1+15t+50t^2+49t^3+15t^4+t^5, \cdots\right)$$

For $m \geq 5$ there are no entries of $C_n^{(m)}(t)$ in OEIS. But

$$\left(c_{n+1}^{(5)}(t)\right) = \left(1, 1+t, 1+t+t^2, 1+2t+2t^2+t^3, 1+2t+4t^2+2t^3+t^4, 1+3t+6t^2+5t^3+3t^4+t^5, \cdots\right),$$

is A276696.

Here we have the recurrence

$$\binom{n}{2k+1}^{(5)} = \binom{n-1}{2k}^{(5)} + \binom{n-2}{2k+1}^{(5)}$$

$$\binom{n}{2k}^{(5)} = \binom{n-1}{2k-1}^{(5)} + \binom{n-1}{2k}^{(5)}.$$

Here we also have $c^{(5)}(t,z)c^{(5)}(-t,-z) = C^{(5)}\left(t^2, z^2\right)$.

For $m = 6$ we get

$$\left(C_{n+1}^{(6)}(t)\right) = \begin{pmatrix} 1, 1+t, 1+3t+t^2, 1+6t+6t^2+t^3, 1+10t+20t^2+10t^3+t^4, \\ 1+15t+50t^2+50t^3+15t^4+t^5, 1+21t+105t^2+174t^3+105t^4+21t^5+t^6, \cdots \end{pmatrix}$$

and

$$\left(c_{n+1}^{(6)}(t)\right) = \left(1, 1+t, 1+t+t^2, 1+2t+2t^2+t^3, 1+2t+4t^2+2t^3+t^4, 1+3t+6t^2+6t^3+3t^4+t^5, \cdots\right),$$

which is Losanitsch's triangle $\left(\binom{n}{k}^{(6)}\right)$, OEIS A034851, which has been studied in [3].

Here we get

$$C^{(6)}_{n+1}(t)=\sum_{2j\le n}\binom{n}{2j}C^{(2)}_j t^j(1+t)^{n-2j},$$

$$c^{(6)}_{n+1}(t)=\sum_{2j\le n}(-1)^j\binom{\left\lfloor \frac{n}{2}\right\rfloor}{j}C^{(2)}_j t^j(1+t)^{n-2j}$$

For $m\ge 2k+1$ we get $\begin{bmatrix} n\\ k\end{bmatrix}^{(m)}=N_{n+1,k}=\binom{n}{k}\binom{n+1}{k}\frac{1}{k+1}$

and $\binom{n}{k}^{(m)}=\left\langle \begin{matrix} n+1\\ k\end{matrix}\right\rangle.$

I could not find closed formulas in the general case. Computations suggest

Con**jecture 2**

For $m\ge 1$

$$\begin{aligned} C^{(4m)}\left(t^2,z^2\right)&=c^{(4m)}\left(t,z\right)c^{(4m)}\left(-t,-z\right),\\ C^{(4m+1)}\left(t^2,z^2\right)&=c^{(4m+1)}\left(t,z\right)c^{(4m+1)}\left(-t,-z\right). \end{aligned} \tag{79}$$

These can be interpreted as analogues of (6).

As an analogue of the fact that the coefficients of $C_{n+1}(t)$ and $c_{n+1}(t)$ are unimodal and palindromic and satisfy

$$C_{n+1}(t)=\sum_{k=0}^{\left\lfloor \frac{n}{2}\right\rfloor}\binom{n}{2k}C_k t^k(1+t)^{n-2k} \tag{80}$$

and

$$c_{n+1}(t)=\sum_{j=0}^{\left\lfloor \frac{n}{2}\right\rfloor}(-1)^j\binom{\left\lfloor \frac{n}{2}\right\rfloor}{j}C_j t^j(1+t)^{n-2j}. \tag{81}$$

computations suggest that the triangles $\left(\begin{bmatrix} n\\ k\end{bmatrix}^{(2m)}\right)$ and $\left(\binom{n}{k}^{(2m)}\right)$ are unimodal and palindromic and satisfy

**Conjecture 3**

$$C_{n+1}^{(2m)}(t)=\sum_{j=0}^{\left\lfloor\frac{n}{2}\right\rfloor}C_j^{(m-1)}\binom{n}{2j}t^j(1+t)^{n-2j},\quad c_{n+1}^{(2m)}(t)=\sum_{j=0}^{\left\lfloor\frac{n}{2}\right\rfloor}(-1)^jC_j^{(m-1)}\binom{\left\lfloor\frac{n}{2}\right\rfloor}{j}t^j(1+t)^{n-2j}.$$

---

Email: johann.cigler@univie.ac.at;  www: johann-cigler.com